\documentclass[12pt]{article}
\usepackage[DIV12]{typearea} 
\usepackage[utf8]{inputenc}
\usepackage{amsmath}
\usepackage{enumerate}
\usepackage{amssymb}
\usepackage{amsfonts}
\usepackage{textcomp}
\usepackage{stmaryrd}
\usepackage[all]{xy}
\DeclareMathOperator*{\coker}{coker}
\DeclareMathOperator*{\rk}{rk}
\DeclareMathOperator*{\Pic}{Pic}

\DeclareMathOperator{\Hom}{Hom}

\DeclareMathOperator{\ch}{ch}

\newcommand{\ta}[1]{#1^{[2]}}

\newcommand{\mc}[1]{\mathcal{#1}}
\newcommand{\tamc}[1]{\ta{\mc{#1}}}

\newcommand{\pf}{\em Proof: \em}
\newcommand{\pfend}{\hfill$\Box$}

\newcommand{\ZZ}{\mathbb{Z}}

\newcommand{\PP}{\mathbb{P}}
\newcommand{\CC}{\mathbb{C}}
\DeclareMathOperator*{\im}{im}

\begin{document}
\newtheorem{thm}{Theorem}[section]
\newtheorem{prop}[thm]{Proposition}
\newtheorem{lem}[thm]{Lemma}
\newtheorem{defi}[thm]{Definition}
\newtheorem{defilemma}[thm]{Definition/Lemma}
\newtheorem{cor}[thm]{Corollary}
\newtheorem{exa}[thm]{Example}
\title{Stability of tautological bundles on the Hilbert scheme of two points on a surface}
\author{Malte Wandel\\\\Leibniz Universität Hannover\\e-mail: wandel@math.uni-hannover.de}

\maketitle
\begin{abstract}
Let $(X,H)$ be a polarised smooth projective surface satisfying $H^1(X,\mc{O}_X)=0$ and let $\mc{F}$ be either a rank one torsion-free sheaf or a rank two $\mu_H$-stable vector bundle on $X$. Assume that $c_1(\mc{F})\neq0$. In this article it is shown that the rank two, respectively rank four tautological sheaf $\tamc{F}$ associated with $\mc{F}$ on the Hilbert square $\ta{X}$ is $\mu$-stable with respect to a certain polarisation.
MCS:\ 14D20, 14J28, 14J60, 14F05
\end{abstract}
\tableofcontents
\setcounter{section}{-1}
\section{Introduction}
Let $X$ be an algebraic K3 surface with polarisation $H\in\Pic X$ and let $v=(r,c,s)\in \ \mathbb{N} \oplus\text{NS}(X)\oplus \mathbb{Z}$. Mukai has shown that in many cases --- if $v$ is carefully chosen --- the moduli space of $H$-semistable sheaves of rank $r$, first Chern class $c$ and second Chern class $s$ is again a smooth compact complex manifold carrying a holomorphic symplectic structure. In fact, all these moduli spaces are deformation equivalent to $\text{Hilb}^n(X)$ for some $n\geq0$. Now the natural question arises what happens if we start with another hyperkähler manifold and study the geometry of moduli spaces of sheaves on this manifold. Not much is known about this topic and one of the fundamental questions is the following: does there exist a symplectic structure on these moduli spaces? Of course answering this question in general will be very complicated. But one could hope for at least finding an example of a such a moduli space that does carry such a symplectic structure. Therefore we need examples of vector bundles on higher dimensional hyperkähler manifolds and then we have to inquire about the stability of these bundles. One big class of examples are the so-called 'tautological bundles' on the Hilbert schemes of K3 surfaces. They arise as the images of vector bundles on a K3 under a Fourier$-$Mukai transform. We will concentrate on the case of $\text{Hilb}^2(X)$, where $X$ is a projective K3-surface. Schlickewei has shown in \cite{Schl} that in many cases tautological bundles associated with line bundles on $X$ are stable with respect to a carefully chosen polarisation on the Hilbert scheme. We will extend this result by showing that, in fact, every rank two tautological sheaf associated with any rank one torsion-free sheaf having non-vanishing first Chern class is stable with respect to some polarisation. Furthermore we will prove that the rank four tautological vector bundle associated with any stable rank two bundle is stable. Again we assume that the first Chern class is nontrivial. This provides us with quite a big variety of stable vector bundles on $\text{Hilb}^2(X)$. In a forthcoming paper it will be shown that in some cases the component of the moduli space of sheaves on $\text{Hilb}^2(X)$ containing the tautological sheaves is smooth and isomorphic to the moduli space of sheaves on the K3 surface. It is therefore an irreducible holomorphic symplectic manifold.\\
In fact, all results concerning the stability of the tautological sheaves are valid for any smooth projective surface $X$ satisfying $h^1(X,\mathcal{O}_X)=0$, so they will be presented in this generality.\\
\\
\em Notations and Conventions: \em
\begin{itemize}
\item In this article all schemes and varieties will be defined over the field of complex numbers.
\item For a vector bundle $\mc{E}$ we write $\PP(\mc{E})$ for $Proj(Sym(\mc{E}^\vee))$ following the definition in Fulton's textbook. In this definition $\PP(\mc{E})$ is the bundle of lines of $\mc{E}$.
\item By $A^\star(Y)$ we denote the Chow ring of any smooth projective variety $Y$.
\end{itemize}
\noindent \em Acknowledgements: \em  I want to thank Klaus Hulek, David Ploog, Marc Nieper-Wißkirchen and Andreas Krug for many useful comments and suggestions. Special thanks go to the unknown referee for helping me improving this paper considerably.

\section{The geometric set-up}\label{sectiongeometricsetup}
Let $X$ be a projective surface satisfying $h^1(X,\mc{O}_X)=0$ and choose a polarisation $H$. Throughout this text we will consider the following basic blow-up and projections diagram:
\[\xymatrix{D \ar[rr]^(.4)i \ar[d]^{\sigma_D} &  & \widetilde{X\times X} \ar[rrd]^\pi \ar[d]^\sigma \ar @/_/ [ldd]_(.4){r_1} \ar  @/^/ [rdd]^{r_2} & & &\ta{X}\times X \ar[ld]^p \ar[rd]^q \\
X \ar[rr]^(.4)\Delta & & X\times X \ar[ld]^{\pi_1} \ar[rd]_{\pi_2} & & \ta{X}&&X \\
  & X & & X &&&.}\]
Here $\Delta$ is the diagonal embedding, $\sigma$ is the blowing-up morphism (we are blowing up the diagonal), $D\simeq \mathbb{P}(\mc{N}_{X|X\times X})\simeq \mathbb{P}(\mc{T}_X)$ denotes the exceptional divisor together with the projection $\sigma_D$, the inclusion $i$ and $\mc{O}_D(1)$, the dual of the tautological line bundle. It is well know that $\mc{N}_{D|\widetilde{X\times X}}\cong \mc{O}_D(-1)$ (see for example Theorem II 8.24 in \cite{Har}). Furthermore $\pi_1$, $\pi_2$, $p$ and $q$ denote the natural projections onto the particular factors and $r_1$ and $r_2$ are the compositions of $\pi_1$ and $\pi_2$ with $\sigma$. Last but not least we have the flat two-to-one covering $\pi$.\\
\\
We will continue with some considerations concerning the Picard groups of the varieties we are looking at. Note that by the assumption $h^1(X,\mc{O}_X)=0$ a line bundle is uniquely determined by its first Chern class. For the same reason this holds true for $X\times X$ and $\widetilde{X\times X}$. We will therefore often use the same notation for a line bunlde as for the corresponding classes in the Chow ring and cohomology.\\
We have $\Pic(X\times X)\cong(\Pic X)^{\boxplus 2}$. Here we apply Exercise III 12.6b) in \cite{Har} since $h^1(X,\mc{O}_X)=0$. And accordingly we have $\Pic(\widetilde{X\times X})\cong(\Pic X)^{\boxplus 2} \oplus \ZZ D$. We will write an element of $\Pic(\widetilde{X\times X})$ as $g\otimes 1+1\otimes h + aD$ for some $g,h \in \Pic X$ and $a\in\ZZ$ and denote the corresponding line bundle by $\mc{L}_{(g,h,a)}$.\\
Furthermore it is well known that $\Pic\ta{X}\cong \Pic X\oplus\ZZ\delta$, where $\delta$ is a class such that $2\delta$ is the exceptional divisor in $\ta{X}$ coming from the blow-up of the diagonal in the quotient $(X\times X)/S_2$. We will denote the line bundle corresponding to $\delta$ by $\mc{L}_\delta$. Thus we can write every element in $\Pic(\ta{X})$ as $\mc{L}_X\otimes \mc{L}_\delta^{\otimes a}$ for some $\mc{L}_X\in\Pic X$ and $a\in\ZZ$. Note that with this notation we have $\pi^\star\mc{L}_X=\mc{L}_{(l,l,0)}$ for a line bundle $\mc{L}_X$ on $X$ with first Chern class $l$. Futhermore we have the relations $\pi^\star\delta=D$ and $\pi_\star D=2\delta$.\\
\\
Next let us summarise the most important facts about the Chow rings of the varieties involved in the upper diagram. We will follow very closely \cite{Ful}, Sections 6.7 and 15.4, especially Lemma 15.4. On $D=\mathbb{P}(\mc{T}_X)$ we have the short exact sequence: $0\rightarrow \mc{O}_D(-1)\rightarrow \sigma_D^\star\mc{N}_{X|X\times X}\rightarrow \mc{Q}\rightarrow 0,$ where $\mc{Q}$ is the universal quotient line bundle. We have $\mc{N}_{X|X\times X}\simeq\mc{T}_X$ and --- by comparing Chern classes --- we can therefore see that $\mc{Q}\simeq\mc{O}_D(1)\otimes \sigma_D^\star\omega_X^\vee$:
\begin{equation}\label{univquot}
0\rightarrow \mc{O}_D(-1)\rightarrow \sigma_D^\star\mc{T}_X\rightarrow \mc{O}_D(1)\otimes \sigma_D^\star\omega_X^\vee\rightarrow 0.
\end{equation}
Let $\xi$ denote the first Chern class of $\mc{O}_D(1)$. By Remark 3.2.4 and Theorem 3.3 in \cite{Ful} we have
\begin{equation*}
A^\ast(D)\cong A^\ast(X)[\xi]/(\xi^2+c_1(\mc{T}_X)\xi +c_2(\mc{T}_X)).
\end{equation*}
Proposition 6.7e) in \cite{Ful} describes the structure of $A^\ast(\widetilde{X\times X})$. We gather the most important identities in this ring in the following lemma.
\begin{lem} \label{lemident} Let $\alpha,\beta,\gamma \in A^\ast(X)$. In $A^\ast(\widetilde{X\times X})$ we have the following identities:
\begin{enumerate}[a)]
\em \item \em
\[i_\star(\xi\cdot\sigma_D^\star(\alpha))=\sigma^\star\Delta_\star(\alpha) +i_\star\sigma_D^\star(\alpha\cdot\omega_X),\]
\em \item \em
\[i^\star i_\star \lambda = -\xi\cdot \lambda\text{, for all }\lambda\in A^\ast(D),\]
\em \item \em 
\[i_\star\sigma_D^\star \alpha\cdot\sigma^\star(\beta\otimes \gamma)=i_\star\sigma_D^\star(\alpha\cdot \beta\cdot \gamma),\]
\em \item \em
\[i_\star\sigma_D^\star(\alpha)\cdot i_\star\sigma_D^\star(\beta)=-\sigma^\star\Delta_\star(\alpha\cdot \beta)-i_\star\sigma_D^\star(\alpha\cdot\beta\cdot\omega_X).\]
\end{enumerate}
\end{lem}
\pf a) Follows from the general formula in Prop. 6.7. a) in \cite{Ful}. Note that in this case the excess normal bundle is just the universal quotient bundle denoted by $\mc{Q}$ above. We have $c_1(\mc{Q})=c_1(\mc{O}_D(1)\otimes \sigma_D^\star\omega_X^\vee)=\xi-\sigma_D^\star\omega_X$.\\
b) This is the self-intersection formula Cor 6.3 in \cite{Ful}:
\[i^\star i_\star\lambda= c_1(\mc{N}_{D|\widetilde{X\times X}})\cdot\lambda=c_1(\mc{O}_D(-1))\cdot\lambda=-\xi\cdot\lambda.\]
c) We have $\alpha\cdot \beta \cdot \gamma = \alpha \cdot \Delta^\star(\beta\otimes \gamma)$. Applying $\sigma_D^\star$ we get
\begin{equation*}
\begin{array}{rrcl}
&\sigma_D^\star (\alpha\cdot \beta \cdot \gamma)& =& \sigma_D^\star(\alpha\cdot \Delta^\star(\beta\otimes \gamma))\\
=&\sigma_D^\star \alpha\cdot \sigma_D^\star \Delta^\star(\beta\otimes \gamma)&=&\sigma_D^\star \alpha\cdot i^\star\sigma^\star(\beta\otimes \gamma).
\end{array}
\end{equation*}
Now we apply $i_\star$ and use the projection formula.\\
d) We use the projection formula and then b) to find
\[
 i_\star\sigma_D^\star(\alpha)\cdot i_\star\sigma_D^\star(\beta)=i_\star(i^\star i_\star\sigma_D^\star(\alpha)\cdot \sigma_D^\star(\beta)) = - i_\star(\xi\cdot\sigma_D^\star(\alpha)\cdot \sigma_D^\star(\beta))=- i_\star(\xi\cdot\sigma_D^\star(\alpha\cdot\beta)).
\]
Now we apply a) and we are done.\pfend\\
\\
\begin{cor}\label{corident} We have
\begin{enumerate}[a)]
\em \item \em
\[i^\star D=-\xi,\]
where we denote $i_\star[D]\in A^3(\widetilde{X\times X})$ simply by $D$,
\em \item \em
\[D^2=-i_\star\xi=-\sigma^\star\Delta-i_\star\sigma_D^\star(\omega_X),\]
 where $\Delta$ also denotes the cohomology class of the diagonal in $X\times X$, and finally
\em \item \em
\[(\sigma^\star\Delta)^2=\sigma^\star\Delta_\star(c_2(\mc{T}_X)).\]
\end{enumerate}
\end{cor}
\pf a) Apply b) of the lemma to $\lambda=[D]$.\\
b) We use a) and for the second equality we apply a) of the Lemma to $\alpha=[X]$ to get
\[D^2=i_\star i^\star D=i_\star(-\xi)=-\sigma^\star\Delta-i_\star\sigma_D^\star(\omega_X).\]
c) Very similarly to the proof of b) in the lemma we use the self-intersection formula:
\[
(\sigma^\star\Delta)^2=\sigma^\star(\Delta^2)=\sigma^\star\Delta_\star\Delta^\star\Delta=\sigma^\star\Delta_\star\Delta^\star\Delta_\star[X]=\sigma^\star\Delta_\star(c_2(\mc{N}_{X|X\times X}))=\sigma^\star\Delta_\star(c_2(\mc{T}_X)).
\]
\pfend\\
\\
Let us finish this chapter by determining the canonical line bundles of $D$ and $\widetilde{X\times X}$. On $\widetilde{X\times X}$ we have a short exact sequence:
\[0\rightarrow \mc{T}_{\widetilde{X\times X}} \rightarrow \sigma^\star\mc{T}_{X\times X}\rightarrow i_\star(\mc{O}_D(1)\otimes \sigma_D^\star\omega_X^\vee)\rightarrow 0.\]
We immediately see $c_1(\mc{T}_{\widetilde{X\times X}})=r_1^\star c_1(\mc{T}_X) + r_2^\star c_1(\mc{T}_X)-D$ and therefore $\omega_{\widetilde{X\times X}}=\mc{L}_{(\omega_X,\omega_X,1)}$.\\
Next, on $D$ we have the exact sequence:
\[0\rightarrow \mc{T}_D \rightarrow i^\star\mc{T}_{\widetilde{X\times X}}\rightarrow \mc{O}_D(-1)\rightarrow 0.\]
Again, we derive $c_1(\mc{T}_D)=2\sigma_D^\star c_1(\mc{T}_X)+2\xi$, so $\omega_D\simeq \sigma_D^\star(\omega_X^\vee)^{\otimes 2}\otimes\mc{O}_D(-2)$.

\section{Tautological bundles}\label{sectiontautbundles}
Now let $\mc{F}$ be a vector bundle on $X$ of rank $r$ with first Chern class $f$. Recall that in $\ta{X}\times X$ there is the \em universal subscheme \em $\Xi$ consisting of pairs $(\xi,x)$ such that $x\in\xi$. We define the \em tautological bundle associated with $\mc{F}$ \em to be the image of $\mc{F}$ under the Fourier$-$Mukai transform with the structure sheaf of the universal subscheme as kernel:
\[ \tamc{F}:= Rp_\star(q^\star\mc{F}\otimes\mc{O}_\Xi).\]
Since we are only considering the case of the second Hilbert scheme we can simplify this definition. Indeed, the universal subscheme $\Xi$ is isomorphic to the blow-up $\widetilde{X\times X}$ of $X\times X$ along the diagonal. A detailed discussion of this fact can be found in \cite[Sect. 1]{EGL}. Via this isomorphism $p$ restricted to $\Xi$ is corresponds to the two-to-one cover $\pi$ and $q$ to the morphism $r_1=\sigma\circ\pi_1$. Thus we end up with the much simpler formula:
\[\tamc{F}=\pi_\star r_1^\star\mc{F}.\]
\em Remark: \em We see immediately that this process is, in fact, an exact functor and we do not need to derive the pushforward along the finite morphism $\pi$.\\\\
Now $\tamc{F}$ is, of course, a vector bundle on $\ta{X}$ of rank $2r$ and we have the following formula for its dual:
\begin{lem}\label{lemdual}
Let $\mc{F}$ be a vector bundle on $X$. Then
\begin{equation}\label{eqtautdual}
(\tamc{F})^\vee \simeq \mc{F}^{\vee[2]}\otimes \mc{L}_\delta.
\end{equation}
\end{lem}
\em Proof: \em Using Grothendieck$-$Verdier duality we have
\begin{eqnarray*}
\begin{array}{rclcl}
\tamc{F}\vee &=& \mc{H}om_{\mc{O}_{\ta{X}}}(\pi_\star r_1^\star\mc{F},\mc{O}_{\ta{X}})  &\simeq &
\pi_\star\mc{H}om_{\mc{O}_{\widetilde{X\times X}}}(r_1^\star\mc{F},\mc{L}_{(\omega_X,\omega_X,1)}\otimes\pi^\star\omega_{\ta{X}}^\vee)\\
&\simeq &
\pi_\star\mc{H}om_{\mc{O}_{\widetilde{X\times X}}}( r_1^\star\mc{F},\mc{L}_{(0,0,1)})
& \simeq & \pi_\star(r_1^\star\mc{F}^\vee\otimes \pi^\star\mc{L}_\delta)\\
&\simeq&\mc{F}^{\vee[2]}\otimes \mc{L}_\delta.
\end{array}
\end{eqnarray*}
Note that we have used here that in the identification $\Pic(\ta{X})\cong\Pic(X)\oplus\ZZ\delta$ we have $\omega_{\ta{X}}\simeq \omega_X$.\pfend\\
\\
The pullback $\pi^\star\tamc{F}$ of a tautological sheaf fits into a basic exact sequence as follows:
\begin{equation}\label{basicexseq}
0\rightarrow \pi^\star\tamc{F} \rightarrow r_1^\star\mc{F}\oplus r_2^\star\mc{F} \rightarrow i_\star\sigma_D^\star\mc{F}\rightarrow 0.
\end{equation}
This sequence was already used by Danila (\cite{Dan}) and Schlickewei (\cite{Schl}) to study tautological sheaves and we refer to loc. cit. for a basic proof of its existence. Furthermore note that the exactness of this sequence is a special case of a more general result due to Scala (cf.\ \cite{Sca}). We want to deduce a simple formula for the first Chern class of $\pi^\star\tamc{F}$. We start with a definition.

\begin{defi}
Let $Y$ be a smooth projective variety and let $\mc{E}$ be a sheaf on $Y$ with $\dim \mathrm{Supp}\,\mc{E}\leq d.$ Let $Z_i$ be the irreducible components of the support of $\mc{E}$ of dimension $d$ and denote by $r_i$ the generic rank of $\mc{E}$ on $Z_i.$ We define the \em $d$-cycle associated with $\mc{E}$ \em to be $Z_d(\mc{E}):=\sum_ir_i[Z_i].$
\end{defi}

\begin{prop}\label{lemc1div}
Let $Y$ be a smooth projective variety, $i\colon W\hookrightarrow Y$ a closed subscheme of dimension $d$ and let $\mc{E}$ be a sheaf on $W.$ We have
\[ch(i_\star\mc{E})=Z_m(i_\star\mc{E})+\text{terms of higher codimension}.\]
\end{prop}
\pf This follows from the generalised Grothendieck$-$Riemann$-$Roch theorem as stated in \cite[Sect.\ 18.3]{Ful}:
By Example 18.3.11 in \cite{Ful} we have
\[\tau_Y(i_\star\mc{E})=Z_d(\mc{E})+\text{terms of higher codimension (t.o.h.c)}.\] 
Now we use Theorem 18.3 to proceed:
\[\tau_Y(i_\star\mc{E})=\ch(i_\star\mc{E})\cap\tau_Y(\mc{O}_Y)=(\ch_d(i_\star\mc{E})+\text{t.o.h.c.})\cap([Y]+\text{t.o.h.c.})=\ch_d(i_\star\mc{E})+\text{t.o.h.c.}\]
For the first inequality we used Thm.\ 18.3(2) and for the second Thm.\ 18.3(5). Note that since $Y$ is smooth, $i_\star\mc{E}$ admits a locally free resolution.\pfend

\begin{cor}
We have
\[
c_1(\pi^\star\tamc{F})=r_1^\star f+r_2^\star f - rD.
\]
\end{cor}
\noindent In the sequel we will analyse conditions such that $\tamc{F}$ is stable. As a main ingredient for this to be possible, we will from now on assume that we are given a polarisation $H$ of $X$ and that $\mc{F}$ is $\mu_H$-stable. More precisely for every subsheaf $\mc{E}\subseteq\mc{F}$ of rank $0<r_\mc{E}<r$ we have
\[ \frac{c_1(\mc{E})\cdot H}{r_\mc{E}}<\frac{c_1(\mc{F})\cdot H}{r}.\]
Next we have to fix a polarisation on $\ta{X}$. This is done as follows: for $N\in\mathbb{N}$ we define $H_N:=NH-\delta$. (Recall that we always use the identification $\Pic(\ta{X})\cong\Pic X\oplus\ZZ\delta$.) This divisor is ample for all sufficiently large $N$, say $N\geq N_0$.\\
Now let us assume that there is a destabilising subsheaf $\mc{E}'\subseteq\tamc{F}$. Pulling back both sheaves via $\pi$ we get an inclusion of sheaves on $\widetilde{X\times X}$:
\[ \pi^\star\mc{E}'=:\mc{E}\subseteq \pi^\star\tamc{F}. \]
Since the slope of a vector bundle is just multiplied by two under the finite pullback $\pi^\star$, $\mc{E}$ is also a destabilising subsheaf of $\pi^\star\tamc{F}$ with respect to the polarisation $\widetilde{H_N}=\pi^\star H_N$ of $\widetilde{X\times X}$. Therefore we will, in fact, consider destabilising subbundles of $\pi^\star\tamc{F}$ which come from $\ta{X}$.\\
As a first step towards any considerations about the stability of a vector bundle $\pi^\star\tamc{F}$, we first have to calculate the slope of a sheaf $\mc{E}$ with respect to the given polarisation. It is defined as
\[\mu_{\widetilde{H_N}}(\mc{E}):=\frac{c_1(\mc{E})\widetilde{H_N}^3}{r_{\mc{E}}},\]
considered as a number by integrating against the fundamental class of $\widetilde{X\times X}$. Thus we first calculate the expansion of $\widetilde{H_N}^3:$
\begin{eqnarray*}
\begin{array}{rcl}
\widetilde{H_N}^3&=& (NH\otimes1 + 1\otimes NH -D)^3 \\
&=&(NH\otimes1+1\otimes NH)^3 -3(NH\otimes1+1\otimes NH)^2D+O(N)\\
&=&3N^3(H^2\otimes H+H\otimes H^2)-3N^2(H^2\otimes1+ 2H\otimes H+1\otimes H^2)D+ O(N). \end{array}
\end{eqnarray*}
Now let $\mc{E}$ be a sheaf on $\widetilde{X\times X}$. We write its first Chern class as $c_1(\mc{E})=g\otimes1 + 1\otimes h + aD$, with $g,h\in\Pic X$ and $a\in \ZZ.$ 
\begin{lem}\label{lemmus}Let $\mc{F}$ be a sheaf on $X$ of rank $r$ and let $\mc{E}$ be a sheaf on $\widetilde{X\times X}$ of rank $r_\mc{E}$. We have the following expansions for the slopes of $\mc{E}$ and $\pi^\star\tamc{F}$:
\begin{eqnarray}
\mu_{\widetilde{H_N}}(\mc{E})&=&\frac{1}{r_{\mc{E}}}\big{\{}3H^2(H.(g+h))N^3+12aH^2N^2\big{\}}+O(N),\label{eqtautmue}\\
\mu_{\widetilde{H_N}}(\pi^\star\tamc{F})&=&\frac{3H^2(H.f)}{r}N^3-6H^2N^2+O(N).\label{eqtautmuf}
\end{eqnarray}
\end{lem}
\pf At first note, that formula (\ref{eqtautmuf}) is just the special case of setting $g=h=f$, $a=-r$ and $r_\mc{E}=2r$ in formula (\ref{eqtautmue}). Next from Lemma \ref{lemident} c)  we deduce that $\sigma^\star(A^i(X\times X))\cdot i_\star\sigma_D^\star(A^j(X))=0$ for $i+j>4$. Thus half of the terms in our computation vanish and we are left with
\begin{eqnarray*}
\begin{array}{rcl}
\widetilde{H_N}^3c_1(\mc{E}) &= &3N^3(H^2\otimes H+H\otimes H^2)(g\otimes1 + 1\otimes h)\\
&&-3N^2(H^2\otimes1+ 2H\otimes H+1\otimes H^2)D\cdot aD+O(N).
\end{array}
\end{eqnarray*}
Finally note that by Lemma \ref{lemident} and Corollary \ref{corident} we have 
\begin{eqnarray*}
(H^2\otimes1)D^2&=& -(H^2\otimes1)\sigma^\star\Delta-(H^2\otimes1)i_\star\sigma_D^\star\omega_X\\
&=&-H^2-i_\star\sigma_D^\star(\underbrace{H^2\cdot\omega_X}_{=0})\\
&=&-H^2
\end{eqnarray*}
and similarly for the terms with $H\otimes H$ and $1\otimes H$.\pfend

\section{Destabilising line subbundles of tautological bundles}\label{sectiondestablinesub}
In this section we will show that for $N\geq N_0$ there exist no $H_N$-destabilising line subbundles $\mc{L}'\subseteq\tamc{F}$ in the case $\mc{F}\not\simeq\mc{O}_X$. So assume that $\mc{L}'$ was such a destabilising line subbundle. The pullback $\mc{L}=\pi^\star\mc{L}'$ of such a line bundle is a destabilising line subbundle of the pullback $\pi^\star\tamc{F}$ with respect to $\widetilde{H_N}$. Composing this inclusion with the one from the basic exact sequence (\ref{basicexseq}) we find $\mc{L}\subseteq r_1^\star\mc{F}\oplus r_2^\star\mc{F}$. We will proceed by showing that $\Hom_{\widetilde{X\times X}}(\mc{L},r_i^\star\mc{F})=0$ for $i=1,2$. The situation is completely symmetric, thus we will focus on $\Hom_{\widetilde{X\times X}}(\mc{L},r_1^\star\mc{F})$. We write the first Chern class of $\mc{L}$ as $c_1(\mc{L})=g\otimes1 + 1\otimes h + aD$ with $g=c_1(\mc{G})$ and $h=c_1(\mc{H})$ for some line bundles $\mc{G}$ and $\mc{H}$ on $X$. In fact, since $\mc{L}$ is coming from $\ta{X}$ this class is invariant under the $S_2$-action, that is, $g=h$. But for later use we will proceed in this generality and denote the line bundle class with first Chern class equal to $g\otimes1 + 1\otimes h + aD$ simply by $\mc{L}_{(g,h,a)}$. We have the following central result:
\begin{prop}\label{prophomef} For all $g,h\in\Pic X$ and $a\in \ZZ$ we have
\begin{equation}
\Hom_{\widetilde{X\times X}}(\mc{L}_{(g,h,a)},r_1^\star\mc{F})\subseteq
\Hom_X(\mc{G},\mc{F})^{h^2(X,\mc{H}\otimes\omega_X)}.
\end{equation}
\end{prop}
\pf Consider the defining exact sequence of the structure sheaf of the exceptional divisor $D$:
\[0\rightarrow \mc{L}_{(0,0,-1)}\rightarrow \mc{O}_{\widetilde{X\times X}}\rightarrow \mc{O}_D\rightarrow 0.\]
Tensoring this sequence with $\mc{L}_{(0,0,a)}$ we have
\begin{eqnarray}
0\rightarrow \mc{L}_{(0,0,a-1)}\rightarrow \mc{L}_{(0,0,a)}\rightarrow \mc{O}_D(-a)\rightarrow 0.\label{definingD}
\end{eqnarray}
So we see immediately that $\sigma_\star\mc{L}_{(0,0,a)}$ is contained in $\mc{O}_{X\times X}$ for all $a\in \ZZ$. Thus we find that $r_{1\star}(r_2^\star\mc{H}^\vee\otimes\mc{L}_{(0,0,-a)})\simeq\pi_{1\star}(\pi_2^\star\mc{H}^\vee\otimes \sigma_\star\mc{L}_{(0,0,-a)})$ is a subsheaf of $\pi_{1\star}\pi_2^\star\mc{H}^\vee\simeq H^0(\mc{H}^\vee)\otimes \mc{O}_X\simeq \mc{O}_X^{h^2(X,\mc{H}\otimes\omega_X)}.$ Now using the projection formula and adjunction we get:
\[
\begin{array}{rcl}
\Hom_{\widetilde{X\times X}}(\mc{L}_{(g,h,a)},r_1^\star\mc{F}) &\cong&
\Hom_{\widetilde{X\times X}}(r_1^\star\mc{G},r_1^\star\mc{F}\otimes \mc{L}_{(0,-h,-a)})\cong\\
\Hom_X(\mc{G},r_{1\star}(r_1^\star\mc{F}\otimes \mc{L}_{(0,-h,-a)}))&\cong&
\Hom_X(\mc{G},\mc{F}\otimes r_{1\star}(r_2^\star\mc{H}^\vee\otimes \mc{L}_{(0,0,-a)})).
\end{array}
\]
 Together with the inclusion of above we are done.\pfend\\
\\
\begin{cor}\label{corlinesub}
Let $\mc{F}$ be a $\mu_H$-stable vector bundle on $X$ of rank $r$ and first Chern class $c_1(\mc{F})=f$. Then $r_1^\star\mc{F}$ contains no line subbundles $\mc{L}_{(g,h,a)}$ satisfying:
\begin{equation}\label{condnolinesub}
H.(g+h)\geq \frac{H.f}{r},
\end{equation}
except the case $r=1$, $h=0$ and $g=f$.
\end{cor}
\pf So let $\mc{L}_{(g,h,a)}$ be a line subbundle of $r_1^\star\mc{F}$ satisfying the hypothesis of the corollary. We will show that $\Hom_X(\mc{G},\mc{F})^{h^2(X,\mc{H}\otimes\omega_X)}=0$ which yields a contradiction to Proposition \ref{prophomef}.\\
If $H.h>0$ we have $0=h^0(X,\mc{H}^\vee)=h^2(X,\mc{H}\otimes\omega_X)$ and we are done.\\
If $H.h\leq0$ we see
\begin{equation}\label{eqhghf}
H.g\geq H.(g+h)\geq\frac{H.f}{r}.
\end{equation}
So if $\mc{G}\not\simeq\mc{F}$ by the stability of $\mc{F}$ we have $\Hom_X(\mc{G},\mc{F})=0$.\\
If $\mc{G}\simeq\mc{F}$ we must have $r=1$ and equalities everywhere in equation (\ref{eqhghf}), so $H.h=0$. But then again $h^2(X,\mc{H}\otimes\omega_X)=0$ for all such $\mc{H}$ but the trivial line bundle, i.e. $h=0$.\pfend\\
\\
Now that we have an explicit description for possible homomorphisms from a line bundle
to $\pi^\star\tamc{F}$, let us have a closer look at the destabilzing condition for line subbundles in $\pi^\star\tamc{F}$.
\begin{lem}\label{lemdestcond}
For sufficiently large $N$ a line subbundle $\mc{L}_{(g,h,a)}$ in $\pi^\star\tamc{F}$ is $H_N$-destabilising if
\begin{eqnarray*}
rH.(g+h)>H.f & \text{or} & rH.(g+h)=H.f \text{ and } a\geq0.
\end{eqnarray*}
\end{lem}
\pf Equation (\ref{eqtautmue}) of Lemma \ref{lemmus} computes the expansion of the slope of $\mc{L}_{(g,h,a)}$ as
\[\mu_{\widetilde{H_N}}(\mc{L}_{(g,h,a)})=3H^2(H.(g+h))N^3+12aH^2N^2 + O(N).\]
And we also derived the expansion of $\mu_{\widetilde{H_N}}(\pi^\star\tamc{F})$ in (\ref{eqtautmuf}) of Lemma \ref{lemmus}:
\begin{eqnarray*}
\mu_{\widetilde{H_N}}(\pi^\star\tamc{F})=\frac{3H^2(H.f)}{r}N^3-6H^2N^2+O(N).
\end{eqnarray*}
Thus $\mc{L}_{(g,h,a)}$ is destabilising if either $rH.(g+h)>H.f$ or $rH.(g+h)=H.f$ and $2a>-1$. Since $a\in\ZZ$ we can replace the last inequality by $a\geq0$.\pfend\\

\begin{thm} \label{thmnodestablinesub}
Let $\mc{F}$ be a $\mu_H$-stable vector bundle on $X$ of rank $r$ and first Chern class $c_1(\mc{F})=f$. Assume $\mc{F}\not\simeq\mc{O}_X$. Then for sufficiently large $N$ the tautological vector bundle $\tamc{F}$ on $\ta{X}$ has no $\mu_{H_N}$-destabilising line subbundles.
\end{thm}
\pf As explained before we are reduced to considering an $S_2$-equivariant destabilising line subbundle $\mc{L}_{(g,g,a)}$ of $\pi^\star\tamc{F}$. The destabilising condition yields $H.g\geq\frac{H.f}{2r}.$ So by Corollary \ref{corlinesub} such a line subbundle cannot exist.\pfend

\section{The cases $r=1$ and $r=2$}\label{sectioncasesr1r2}
From Theorem \ref{thmnodestablinesub} we deduce:
\begin{cor}\label{corstablb}
Let $\mc{F}$ be a line bundle on $X$ not isomorphic to $\mc{O}_X$. Then for sufficiently large $N$, $\tamc{F}$ is a $\mu_{H_N}$-stable rank two vector bundle on $\ta{X}$.
\end{cor}
\pf Since $\tamc{F}$ has rank two we only have to consider torsion-free destabilising subsheaves of rank one. If $\mc{E}$ is such a subsheaf we can embed it into its reflexive hull $\mc{E}^{\vee\vee}$. This is a reflexive rank one sheaf, i.e. a line bundle. Since $\tamc{F}$ is locally free it is also reflexive. Now $\mc{E}^{\vee\vee}$ is a subbundle of $\tamc{F}$ and the first Chern classes of $\mc{E}^{\vee\vee}$ and $\mc{E}$ coincide. Therefore $\mc{E}^{\vee\vee}$ is destabilising. This gives a contradiction to Theorem \ref{thmnodestablinesub}.\pfend\\
\\We can generalise this result to arbitrary torsion free rank one sheaves on $X$ with nonvanishing first Chern class:
\begin{thm}
Let $\mc{F}$ be a torsion free rank one sheaf on $X$ satisfying $c_1(\mc{F})\neq0$. Then for sufficiently large $N$, $\tamc{F}$ is a $\mu_{H_N}$-stable rank two torsion free sheaf on $\ta{X}$.
\end{thm}
\pf Every torsion free rank one sheaf $\mc{F}$ on a surface can be written as $\mc{F}\simeq\mc{L}\otimes\mc{I}_Z$ for some line bundle $\mc{L}$ and an ideal sheaf $\mc{I}_Z$ of a zero dimensional subscheme $Z\subset X$. We thus have an injection $\mc{F}\subseteq\mc{L}$ and, of course, $c_1(\mc{F})=c_1(\mc{L})$. In particular, the line bundle $\mc{L}$ is not trivial.\\
Now since $(-)^{[2]}$ is an exact functor we see (cf.\ Lemma 23 in \cite{Sca}) that $\tamc{F}$ is also torsion free. Furthermore we have an injection $\tamc{F}\subseteq\tamc{L}$. But $c_1(\tamc{F})=c_1(\tamc{L})$ because the cokernel of the inclusion $\tamc{F}\hookrightarrow\tamc{L}$ is $\tamc{O}_Z$ which is supported in codimension two. So the stability of $\tamc{F}$ follows immediately. \pfend\\
\\
Now we want to consider the case $r=\rk\mc{F}=2$. We have seen before that $\tamc{F}$ cannot contain destabilising line subbundles. In this section we will prove that in most cases, in fact, $\tamc{F}$ does not contain any destabilising subsheaves. We start with a technical lemma.
\begin{lem}\label{lemmaxsat}
Let $(Y,\mc{O}(1))$ be a polarised smooth projective variety and let $\mc{H}$ be a pure sheaf on $Y$. The maximal destabilising subsheaf of $\mc{H}$ is saturated.
\end{lem}
\pf Denote by $\mc{H}'$ the maximal destabilising subsheaf. Its saturation $\mc{H}'_{sat}$ in $\mc{H}$ is a subsheaf of $\mc{H}$ of the same rank containing $\mc{H}'$. Therefore $\mu(\mc{H}'_{sat})\geq\mu(\mc{H}')$. By the maximality of $\mc{H}'$ we must have $\mc{H}'_{sat}\simeq\mc{H}'$.\pfend
\begin{thm}
Let $\mc{F}$ be a rank two $\mu_H$-stable vector bundle on $X$ and assume $f=c_1(\mc{F})\neq0$. Then for sufficiently large $N$, $\tamc{F}$ is a $\mu_{H_N}$-stable rank four vector bundle on $\ta{X}$.
\end{thm}
\pf Let $\mc{E}$ be the maximal destabilising subsheaf of $\pi^\star\tamc{F}$. It is semistable and $S_2$-linearised. Similarly to the proof of Corollary \ref{corstablb} one can show that $\mc{E}$ is reflexive and by Lemma \ref{lemmaxsat} above we see that it is saturated. By Theorem \ref{thmnodestablinesub}, $\mc{E}$ cannot have rank one. So let us first consider the case $\rk\mc{E}=3$ and let us have a look at the corresponding short exact sequence on $\widetilde{X\times X}$:
\[0\rightarrow \mc{E}\rightarrow \pi^\star\tamc{F}\rightarrow \mc{Q} \rightarrow 0,\]
where $\mc{Q}$ is the corresponding destabilising quotient. Let us write $c_1(\mc{E})=e\otimes 1+1\otimes e +aD$.  Using equation (\ref{eqtautdual}) we see that the dual of this sequence looks as follows:
\[0\rightarrow \underbrace{\mc{H}om_{\mc{O}_{\widetilde{X\times X}}}(\mc{Q},\mc{O}_{\widetilde{X\times X}})}_{=:\mc{Q}'}\rightarrow \pi^\star(\mc{F}^{\vee{[2]}})\otimes \mc{L}_{(0,0,1)} \rightarrow \mc{E}^\vee \rightarrow \mc{E}xt^1_{\mc{O}_{\widetilde{X\times X}}}(\mc{Q},\mc{O}_{\widetilde{X\times X}}) \rightarrow 0.\]
Since $\mc{E}$ is saturated, then $\mc{Q}$ is torsion free and so the support of $\mc{E}xt^1_{\mc{O}_{\widetilde{X\times X}}}(\mc{Q},\mc{O}_{\widetilde{X\times X}})$ has codimension at least $2$, so vanishing first Chern class. We compute
\begin{eqnarray*}
c_1(\mc{Q}')&=&c_1\big{(}\pi^\star((\mc{F}^\vee)^{[2]})\big{)}+c_1(\mc{L}_{(0,0,1)})\cdot\rk\big{(}\pi^\star((\mc{F}^\vee)^{[2]})\big{)}-c_1(\mc{E}^\vee)\\
&=&(e-f)\otimes 1 + 1\otimes (e-f)+(a+2)D.
\end{eqnarray*}
Now we may assume that $\mc{Q}'$ is reflexive, i.e. locally free. (If necessary we replace $\mc{Q}'$ by its reflexive hull which still gives a subsheaf of $\pi^\star(\mc{F}^{\vee{[2]}})\otimes \mc{L}_{(0,0,1)}$ with the same first Chern class.) More precisely, $\mc{Q}'\simeq \mc{L}_{(e-f,e-f,a+2)}$. We have an inclusion $\mc{Q}'\otimes\mc{L}_{(0,0,-1)}\hookrightarrow \pi^\star\mc{F}^{\vee[2]}$. Now, by Lemma \ref{lemmus}, the destabilising condition on $\mc{E}$ implies:
\[4H.e\geq3H.f.\]
Thus $2H.(e-f)\geq -\frac{H.f}{2}$ and since $\mc{Q}'\otimes\mc{L}_{(0,0,-1)}\subset r_1^\star\mc{F}$, then by Corollary \ref{corlinesub} we get a contradiction.\\
Finally assume that the maximal destabilising subsheaf of $\pi^\star\tamc{F}$ is a rank two sheaf $\mc{E}$. Again its first Chern class can be written as $c_1(\mc{E})=e\otimes1+1\otimes e + aD$ with $e\in\Pic X$ and $a\in \ZZ$ and by the fundamental exact sequence (\ref{basicexseq}) we get an injective $S_2$-equivariant homomorphism $\mc{E}\hookrightarrow r_1^\star\mc{F}\oplus r_2^\star\mc{F}$. We will denote its composition with the projection onto the first factor by $\beta\colon \mc{E} \rightarrow r_1^\star\mc{F}$. Now we will distinguish three cases:\\
a) rank $\ker\beta=0.$\\
So $\ker\beta$ is a torsion subsheaf of $\mc{E}$, so it is trivial since $\mc{E}$ is torsion-free. So $\beta$ is an isomorphism away from an effective divisor $j\colon Y\hookrightarrow\widetilde{X\times X}$. Thus $\coker\beta$ can be written as $j_\star\mc{K}$ for some sheaf $\mc{K}$ on $Y$. Let $Y=\bigcup_i Y_i$ be the decomposition into irreducible components, then by Proposition \ref{lemc1div} we can write its first Chern class as $c_1(\coker\beta)=\sum_i (Y_i\cdot \rk\mc{K}_i)$, where $\mc{K}_i$ is the restriction of $\mc{K}$ to $Y_i$. On the other hand we can compute the first Chern class of $\coker\beta$ directly:
\[c_1(\coker\beta)=c_1(r_1^\star\mc{F})-c_1(\mc{E})=f\otimes 1-e\otimes1-1\otimes e - aD.\]
Now $Y$ is effective. So if $\rk\mc{K}_i\neq 0$ for some $i$ we must have $(f-e)\otimes1-1\otimes e -aD$ effective. Evaluating against the polarisation $\widetilde{H_N}$ yields $2H.e<H.f$. Together with the destabilising condition on $\mc{E}$ $-$ which implies $2.He\geq H.f$ $-$ we get a contradiction. If $\rk\mc{K}_i=0$ $\forall i$, i.e. $c_1(\coker\beta)=0$ we must have $f=0$ which we excluded.\\
b) rank $\ker\beta=2.$\\
This says that on an open subset $\beta$ has to vanish which by symmetry contradicts the fact that $\mc{E}$ injects into $r_1^\star\mc{F}\oplus r_2^\star\mc{F}$.\\
c) rank $\ker\beta=1.$\\
Now $\im\beta$ is a rank one quotient sheaf of $\mc{E}$ and we write its first Chern class $c_1(\im\beta)=g\otimes1+1\otimes h + bD$. The semistability of $\mc{E}$ yields
\[H.e\leq H.(g+h).\]
At the same time $\im\beta$ is a rank one subsheaf of $r_1^\star\mc{F}$. Denote by $\im\beta^{\vee\vee}$ its reflexive hull. This is a reflexive rank one sheaf, thus a line bundle. And it has the same first Chern class as $\im\beta$, so $\im\beta^{\vee\vee}=\mc{L}_{(g,h,b)}$. The destabilising condition on $\mc{E}$ implies $2H.e\geq H.f$. Putting things together we find a line subbundle $\mc{L}_{(g,h,b)}$ in $r_1^\star\mc{F}$ satisfying $2H.(g+h)\geq H.f$. This is a contradiction to Corollary \ref{corlinesub}.\pfend

\section{The case of the trivial line bundle}\label{sectiontriviallinebundle}
In the previous section we explicitly excluded the case $\mc{F}\simeq \mc{O}_X$. In fact, we have the following result:
\begin{prop}
The tautological vector bundle $\tamc{O}_X$ associated with the trivial line bundle $\mc{O}_X$ is not $\mu_{H_N}$-stable for sufficiently large $N$.
\end{prop}
\pf By \cite[Th\'{e}or\`{e}me 1]{Dan} we have $H^0(\tamc{O}_X)\cong\CC$. Thus the structure sheaf $\mc{O}_{\ta{X}}$ is a line subbundle of $\tamc{O}_X$. We compare the slopes in order to show that $\mc{O}_{\ta{X}}$ is destabilising. By Lemma \ref{lemmus} we have
\begin{eqnarray*}
\mu_{\widetilde{H_N}}(\pi^\star\mc{O}_{\ta{X}})&=&0 \hspace{10pt}\text{ and }\\
\mu_{\widetilde{H_N}}(\pi^\star\tamc{O}_X)&=&-6H^2N^2+O(N).
\end{eqnarray*}
Thus we see that for sufficiently large $N$ the subbundle $\mc{O}_{\ta{X}}$ is destabilising.\pfend

\end{document}